\documentclass[11pt]{article}

\usepackage{amsthm,amssymb,latexsym,graphics,amscd,amsmath,amscd,enumerate,color}

\newcommand{\C}{\mathbb{C}}

\newcommand{\ind}{\mbox{ind}}
\newcommand{\mb}{\mathbb}
\newcommand{\mc}{\mathcal}
\newcommand{\R}{\mathbb{R}}

\newcommand{\Sc}{\mb{S}_{\mb{C}}}
\newcommand{\Sp}{\mbox{Spin}}

\newcommand{\vol}{\mbox{vol}}
\newcommand{\rk}{\mbox{rk}}

\begin{document}

\title{On the eigenvalues of the twisted Dirac operator}
\author{Marcos Jardim \\ IMECC - UNICAMP \\
Department of Mathematics \\ C. P. 6065 \\
13083-970 Campinas-SP, Brazil \\ and \\
Rafael F. Le\~ao \\ Federal University of Paraná \\ Department of Mathematics \\ C. P. 019081 \\ 81531-990 Curitiba-PR, Brazil}

\maketitle

\begin{abstract}
Given a compact Riemannian spin manifold with positive scalar curvature, we find a family of connections $\nabla^{A_t}$ for $t\in[0,1]$ on a trivial vector bundle of sufficiently high rank, such that the first eigenvalue of the twisted Dirac operator $D_{A_t}$ is nonzero and becomes arbitrarily small as $t\to1$. However, if one restricts the class of twisting connections considered, then nonzero lower bounds do exist. We illustrate this fact by establishing a nonzero lower bound for the Dirac operator twisted by Hermitian-Einstein connections over Riemann surfaces.
\end{abstract}


\section{Introduction}
Given a $\Sp$ manifold $(M,g)$, one can consider the spinor bundle $\mb{S}$ and the associated Dirac operator $D$. There is a vast literature about this operator, in particular concerning the behavior of its spectrum. One of the most studied problems is to find lower bounds for the eigenvalues of $D$. The most well-known result was obtained by Friedrich \cite{Fr}, and can be stated as follows: if $(M,g)$ is a compact Riemannian $\Sp$ manifold, with positive scalar curvature $R$, then the eigenvalues of the associated Dirac operator satisfy the inequality
\begin{displaymath}
  \lambda^2 \geq \frac{1}{4} \frac{n}{n-1} R_0
\end{displaymath}
where $R_0$ is the minimum of the scalar curvature and $n$ is the dimension of $M$. If the equality is satisfied, then the scalar curvature $R$ is constant and $M$ is an Einstein manifold. In addition, if further geometric structures in $(M,g)$ are conjecture, the above lower bound can be improved. For example, when $(M,g)$ is a compact K\"ahler manifold, then Kirchberg \cite{Ki}, proved that the eigenvalues of the Dirac operator satisfy the inequality
\begin{displaymath}
  \lambda^2 \geq
  \left\{
  \begin{matrix}
    \frac{1}{4} \frac{k+1}{k} R_0 & \mbox{if $k = \dim_{\C} M$ is odd,} \\
    \/ & \/ \\
    \frac{1}{4} \frac{k}{k-1} R_0 & \mbox{if $k = \dim_{\C} M$ is even,}
  \end{matrix}
  \right.
\end{displaymath}
The influence on the eigenvalues of the Dirac operator of other geometric structures on $M$ are found in \cite{Kr1,Kr2}; see also the survey \cite{JL1}. Kirchberg has also considered estimates in terms of other curvature tensors, see \cite{Ki2} and the references therein.

However, this is not the only Dirac operator one can define on a $\Sp$ manifold $(M,g)$. For instance, consider a Hermitian vector bundle with a compatible connection $(E,\nabla^A)$ over $(M,g)$. Using the connection $\nabla^A$ on $E$, one can define the twisted Dirac operator $D_A$ on $\mb{S} \otimes E$. This operator is extremely important in classical field theory, for it describes particles like electrons and neutrinos coupled to external gauge fields. Very little is known about the behavior of the eigenvalues of $D_A$ in terms of the coupling connection $\nabla^A$.

More precisely, fix the base manifold and geometry $(M,g)$, and a hermitian vector bundle $E\to M$. Let $\mc{A}$ be the affine space of compatible connections on $E$, and consider the functional
\begin{displaymath} \label{functional}
  \begin{split}
    \lambda &: \mc{A} \rightarrow \R^+ \\
    \nabla^A &\mapsto \mid \lambda_1(\nabla^A) \mid
  \end{split}
\end{displaymath}
that associates to each connection $\nabla^A\in\mc{A}$ the absolute value of the first nonzero eigenvalue of the associated Dirac operator $D_A$. In analogy with the untwisted case, it would be interesting to determine whether this functional possesses a nonzero lower bound, and whether it is bounded above.

In \cite{VW}, Vafa and Witten have proved that the functional $\lambda$ does admit an universal upper bound which depends only in the geometry of the compact base manifold $(M,g)$ but not on the twisting bundle $E$; see also \cite{At} for a clear geometrical proof. On the other hand, lower bounds are known only for very special cases, see for instance \cite{Mi}, and rely on strong restrictions both on the base manifold $(M,g)$ and especially on the connection $\nabla^A$.

Despite these particular results, the problem of how the bounds for eigenvalues of the twisted Dirac operator, and consequently the spectrum, changes with the twisting connection is not well understood. The first natural question is whether the functional $\lambda$ has a nonzero lower bound.

In this article, we construct a general example showing that the first eigenvalue of the twisted Dirac operator depends strongly on the twisting connection. More precisely, if $(M,g)$ is a compact Riemannian $\Sp$ manifold with positive scalar curvature $R>0$, then we find a trivial bundle $\underline{\C}^N\to M$, for some $N$ large enough, and a family of connections $\nabla^{A_t}$ on $\underline{\C}^N$ for $t\in[0,1]$, such that the first eigenvalue of the twisted Dirac operator $D_{A_t}$, is nonzero and becomes arbitrarily small as $t\to1$. In other words, a Vafa-Witten type of result for an universal lower bound for the first eigenvalue of the twisted Dirac operator is impossible; lower bounds for the first nonzero eigenvalue necessarily depends on the topology and the geometry of the twisting bundle.

However, if one restricts the class of twisting connections considered, then lower bounds for the nonzero eigenvalues do exist. We illustrate this fact by establishing a lower bound for the Dirac operator twisted by Hermitian-Einstein connections over Riemann surfaces. More generally, we believe that interesting lower bounds for the twisted Dirac operator can be obtained whenever the twisting connection satisfies some classical field equation, like the anti-self-duality Yang-Mills equation.

\paragraph{Acknowledgments.}
M.J. is partially supported by the CNPQ grant number 305464/2007-8 and the FAPESP grant number 2005/04558-0. R.F.L.'s research was supported by a CNPQ doctoral grant.


\section{Arbitrarily Small Eigenvalues}

Let $(M,g)$ be a compact Riemannian $\Sp$ manifold with positive scalar curvature $R$. In this Section, we will construct a 1-parameter family of connections $\nabla^{A_t}$ for $t\in[0,1]$ on a trivial bundle $\underline{\C}^N={\C}^N\times M$ for which the first eigenvalue of the associated twisted Dirac operator $D_{A_t}$, is nonzero and arbitrarily small. The main idea is to show that the trivial bundle $\underline{\C}^N$ admits two connections, $\nabla^0$ and $\nabla^1$, such that the Dirac operator twisted by $\nabla^0$ has trivial kernel while the Dirac operator twisted by $\nabla^1$ has non trivial kernel.

The difficult part is to find the right trivial bundle $\underline{\C}^N$ and the connection $\nabla^1$, after constructing this connection the connection $\nabla^0$ is easily find. To do this, let us first assume that the manifold $(M,g)$ is even dimensional; let $E\to M$ be a hermitian vector bundle over $M$, and $\nabla^A$ be a compatible connection on $E$. In this case, the spinor bundle splits as follows
\begin{equation}
  \mb{S} = \mb{S}^+ \oplus \mb{S}^-~,
\end{equation}
which implies that the twisted Dirac operator also splits $D_A=D_A^+\oplus D_A^-$, where 
\begin{equation}
  \begin{split}
    D_A^+:\mb{S}^+ \otimes E \rightarrow \mb{S}^- \otimes E \\
    D_A^-:\mb{S}^- \otimes E \rightarrow \mb{S}^+ \otimes E
  \end{split}
\end{equation}

The index for the operator $D_A^+$ can be topologically calculated with the following formula
\begin{equation}
  \ind( D_A^+ ) = (-1)^n \int_M ch(E) \wedge \hat{A}(M) \label{ind_aco}
\end{equation}
where $\hat{A}$ is the so-called \^A-genus of $M$, a characteristic class that can be written in terms of the Pontrjagin classes, and $ch(E)$ is the Chern character of the bundle $E$. This expression, together with the analytical index
\begin{equation}
  \ind( D_A^+ ) = \dim \ker D_A^+ - \dim \ker D_A^-
\end{equation}
and the fact that both ways to calculate the index are equal by the Atiyah-Singer index theorem, allows us to conclude that we can find a hermitian vector bundle $E$ with connection $\nabla^A$ such that $D_A$ has nontrivial kernel. Indeed, if the index of $D_A^+$ is not zero, then $D_A$ has a nontrivial kernel. In other words, if the Dirac operator does not have kernel, then the index of $D_A^+$ must vanish.

If $\dim M=2n$, take $E\to M$ to be the pullback of the generating bundle $H$ over $S^{2n}$ by a map $f:M\to S^n$ of degree 1. Then the only non vanishing Chern class of $E$ is the top Chern class $c_n(E)$, which can be represented, in de Rham cohomology, by a generator of $H^{2n}(M)$. Using this construction and the assumption that the original Dirac operator $D$ does not have kernel it is easy to see, using the topological index formula (\ref{ind_aco}), that $D_A^+$ has nonzero index, hence $D_A$ must have kernel for any connection $\nabla^{A}$ we choose in $E$.

Since $M$ is compact, we can find a bundle $E^{\prime}$ such that
\begin{equation}
  E \oplus E^{\prime} \simeq \underline{\C}^N
\end{equation}
for some large enough $N$. This is the desired trivial bundle. Choosing any connection $\nabla^B$ in $E^{\prime}$, set $\tilde{A}=A\oplus B$, which is a connection on $E \oplus E^{\prime} \simeq \underline{\C}^N$. Next, consider the twisted Dirac operator $D_{\tilde{A}} = D_{A}\oplus D_B$; the fact that $D_A$ has nontrivial kernel in $E$ implies that $D_{\tilde{A}}$ also does.

If the dimension of $M$ is odd then we can consider the even dimensional manifold $M \times S^1$ and use the fact that the spectrum of the Dirac of this product manifold is given by
\begin{equation}
  \pm \sqrt{\lambda_j^2 + \beta_k^2}
\end{equation}
where $\lambda_j$ are the eigenvalues of the Dirac operator on $M$ and $\beta_k$ the eigenvalues of the Dirac operator on $S^1$. But the first eigenvalue of the Dirac operator on $S^1$ is zero, so the proof for the odd dimensional case follows as above.

The construction of the connection $\nabla^0$ is trivial. Note that sections of the trivial bundle $\mb{S}\otimes\underline{\C}^N$ can be written in the form
\begin{equation}
  \psi =
  \begin{pmatrix}
    \psi_1 \\
    \psi_2 \\
    \vdots \\
    \psi_N
  \end{pmatrix}
\end{equation}
which implies that, for the trivial connection $\underline{d}$ on $\underline{\C}$, the associated twisted Dirac operator $D_0$ can be written as
\begin{equation}
  D_0 \psi =
  \begin{pmatrix}
    D \psi_1 \\
    D \psi_2 \\
    \vdots \\
    D \psi_N
  \end{pmatrix}
\end{equation}
where $D$ is the free, untwisted Dirac operator associated to $(M,g)$. In particular, this means that the only effect of the coupling with the trivial connection is to change the multiplicity of the eigenvalues, leaving the spectrum unchanged. Since we assumed that the scalar curvature of $(M,g)$ is positive, the Weitzenb\"ock formula
\begin{equation}
  D^2 = \Delta + \frac{1}{4}R
\end{equation}
ensures that the Dirac operator associated to $(M,g)$ does not have kernel and, consequently, that the twisted operator $D_0$ also has trivial kernel.

Summing up, if $(M,g)$ is a Riemannian $\Sp$ manifold with positive scalar curvature, then one can choose a connection, the trivial one, on the trivial vector bundle $\underline{\C}^N$ such that the twisted Dirac operator does not have a kernel.

Now we have in $\underline{\C}^N$ the two connections with the desired properties: one is the trivial connection $\underline{d}$ and the other is $\nabla^{\tilde{A}} = \nabla^{A} \oplus \nabla^B$. The Dirac operator associated to the first has trivial kernel, while the one associated to the latter has nontrivial kernel. Since the space of connections over a fixed vector bundle is an affine space modeled in the space of 1-forms over $M$ with values in the endomorphism bundle, we can write
\begin{equation}
  \nabla^{\tilde A} = \underline{d} + \alpha
\end{equation}
where $\alpha$ is an 1-form over $M$ with values in the endomorphism bundle. This can be used to define the family of connections
\begin{equation}
  \nabla^{A_t} = \underline{d} + t \alpha
\end{equation}
and the family of twisted Dirac operators $D_{A_t}$ in the obvious way. Clearly, $D_{A_0}=D_0$ and $D_{A_1}=D_{\tilde A}$, and by construction, $D_{A_0}$ does not have kernel while $D_{A_1}$ does. Since the functional $\lambda$ defined in (\ref{functional}) is continuous (see \cite{At}), we conclude that $\lambda(D_{A_t})$ becomes arbitrarily small as $t\to 1$, as desired.


\section{Uniform Bound for Riemann surfaces}

In this section we show that for suitable conditions on the base manifold and on the twisting connection $\nabla^A$, it is possible to find lower bounds for the first nonzero eigenvalue of $D_A$.

First of all, note that the topological index formula (\ref{ind_aco}) always implies that the twisted operator $D_A$, over Riemann surfaces, has non vanishing index, which implies that for Riemann surfaces $D_A$ always have non trivial kernel. Because of this fact, lower bounds only makes sense for non null eigenvalues.

Riemann surfaces are naturally K\"ahler manifolds, with the K\"ahler form given by the volume form $\omega = i \xi \wedge \bar{\xi}$. We can use this to make the connection $\nabla^A$ satisfies a compatibility condition known as the Hermitian-Einstein condition. A connection $\nabla^A$ is called a Hermitian-Einstein connection if
\begin{equation}
  \omega \lrcorner F_A = -ic \mb{I}
\end{equation}
where $F_A$ denotes the curvature 2-form of $\nabla^A$ and $\omega \lrcorner$ is the contraction by the K\"ahler form. For background on the importance of the Hermitian-Einstein condition, see \cite{Ko1}.

Another important feature of a Riemann surface $M$ is that the spinor bundle can be explicitly described in terms of forms; more precisely, it is well-known that the spinor bundle associated with the complex structure is $\Sc \simeq \wedge^{(0,0)} \oplus \wedge^{(0,1)}$, so the usual spinor bundle is
\begin{equation}\label{spins}
  \mb{S} \simeq \Sc \otimes K_M^{\frac{1}{2}} \simeq \left( \wedge^{(0,0)} \oplus \wedge^{(0,1)} \right) \otimes k_M^{\frac{1}{2}}
\end{equation}
where $K_M$ is the canonical bundle of $M$. Furthermore, the complex Dirac operator $D$ coincides with a twisted real Dirac $\cal D_S$, where $S$ is the connection on $K_M^{-\frac{1}{2}}$ induced by the Chern connection on $M$.

Now let $E\to M$ be a holomorphic vector bundle of negative degree on $M$ that admits a Hermitian-Einstein connection $\nabla^A$. In \cite{JL2}, the authors have shown that the eigenvalues of the twisted complex Dirac operator satisfy the following lower bound:
\begin{equation}\label{lb-d}
  \lambda^2 \geq - \frac{4 \pi \deg(E)}{\rk(E) \vol(M)}~.
\end{equation}
This lower bound for the eigenvalues of the twisted complex Dirac operator can be translated into a lower bound for the eigenvalues of the twisted real Dirac operator in the following manner. First, consider $L=K_{M}^{\frac{1}{2}}\otimes E$, so that $\deg(L) = \deg(E) - \rk(E)(1-g)$, where $g$ is the genus of $M$. Let $\nabla^{\tilde A}$ be the tensor connection $\nabla^S \otimes \mb{I} + \mb{I} \otimes \nabla^A$. Then the twisted complex Dirac operator $D_{\tilde A}$, on $\Sc \otimes L$, coincides with the twisted real Dirac operator ${\cal D}_A$, on $\mb{S} \otimes E = \Sc \otimes K_{M}^{-\frac{1}{2}}\otimes E$. Thus applying the lower bound (\ref{lb-d}) to the bundle $L$, we obtain
\begin{equation}\label{lb-d2}
  \lambda^2 \geq \frac{4\pi(1-g)}{\vol(M)}- \frac{4 \pi \deg(E)}{\rk(E) \vol(M)}~,
\end{equation}
which is a lower bound for the eigenvalues of the twisted real Dirac operator on $M$. Of course, this lower bound is only significant provided $\deg(E)<\rk(E)(1-g)$.

The Gauss-Bonnet formula can be used to re-write (\ref{lb-d2}) in the following manner:
$$ \lambda^2 \geq \frac{R_0}{2} \left( 1 - \frac{\deg(E)}{(1-g)\rk(E)} \right)~~ {\rm if} ~~ g\ne1 ~,$$
$$ \lambda^2 \geq - \frac{4 \pi \deg(E)}{\rk(E) \vol(M)} ~~ {\rm if} ~~ g=1 ~.$$
Comparing with the results in \cite{Al} (see Theorems 5.3, 5.10 and 5.22 there), we conclude that the above estimates are actually attained in the case where $M$ has constant scalar curvature, $E\to M$ is a line bundle and the connection on $E$ ha constant curvature.


\section{Conclusion}

We have constructed a general example showing that the first eigenvalue of the twisted Dirac operator on a vector bundle over an arbitrary compact Riemannian spin manifold $M$ with positive scalar curvature depends strongly on the twisting connection. This was done by providing a 1-parameter family of connections $\nabla^{A_t}$ on a trivial vector bundle over $M$ of sufficiently high rank such that the first eigenvalue of the twisted Dirac operator $D_{A_t}$, is nonzero and becomes arbitrarily small as $t\to1$. Therefore, one cannot expect absolute lower bounds for the nonzero eigenvalues of a Dirac operator twisted by arbitrary connections.

However, if one imposes conditions on the class of twisting connections, e.g. the connection is satisfy some classical field equation like the anti-self-duality Yang-Mills equation of the Hermitian-Einstein condition, then lower bounds do exist. We illustrate this fact by establishing a sharp, nonzero lower bound for the Dirac operator twisted by Hermitian-Einstein connections over Riemann surfaces.


\end{document}